\documentclass[11pt,fleqn]{article}
\usepackage{amsmath,amssymb,amsthm,mathrsfs}
\usepackage{graphicx}

\newcommand{\eps}{\varepsilon}
\newcommand{\phy}{\varphi}
\newcommand{\R}{\mathbb{R}}
\newcommand{\Z}{\mathbb{Z}}
\newcommand{\T}{\mathbb{T}}
\newcommand{\dee}{\mathrm{d}}
\newcommand{\imag}{\mathrm{i}}
\newcommand{\rem}[1]{}
\newcommand{\X}{{\mathscr{X}}}
\newcommand{\pscal}[2]{\langle #1,#2\rangle}
\renewcommand{\S}{{\cal T}}

\newtheorem{teo}{Theorem}  
\newtheorem{lem}[teo]{Lemma}

\begin{document}

\title{Vanishing Twist near Focus-Focus Points} 
\author{Holger
  R.~Dullin$^{1,2}$, V\~{u} Ng\d{o}c San$^{3}$ \thanks{This research
    was supported by the EPSRC under contract GR/R44911/01 and by the
    European Research Training Network {\it Mechanics and Symmetry in
      Europe\/} (MASIE), HPRN-CT-2000-00113.}  \\ $^{1}$ Fachbereich
  1, Physik, Universit\"at Bremen \\ 28334 Bremen, Germany
\\ $^{2}$ on leave from Department of Mathematical Sciences, 
\\ Loughborough University, LE11 3TU, UK \\ {\small
  H.R.Dullin@lboro.ac.uk} \\ $^{3}$ Institut Fourier, Universit\'e
Grenoble \\ 38402 Saint Martin d'H\`eres, France \\ {\small
  svungoc@ujf-grenoble.fr} }

\maketitle

\begin{abstract}
  
  We show that near a focus-focus point in a Liouville integrable
  Hamiltonian system with two degrees of freedom lines of locally
  constant rotation number in the image of the energy-momentum map are
  spirals determined by the eigenvalue of the equilibrium.  From this
  representation of the rotation number we derive that the twist
  condition for the isoenergetic KAM condition vanishes on a curve in
  the image of the energy-momentum map that is transversal to the line
  of constant energy. In contrast to this we also show that the
  frequency map is non-degenerate for {\em every} point in a
  neighborhood of a focus-focus point.

\vspace*{1ex}

\noindent
Math. Class.: 37J35, 37J15,  37J40, 70H06, 70H08, 37G20

  

\end{abstract}

\newpage

\section{Introduction}

Consider a neighborhood of an equilibrium point of a Hamiltonian
system with two degrees of freedom.  In local symplectic coordinates
$(x,y,p_x,p_y) \in \R^4$ with standard symplectic structure $\Omega =
\dee x\wedge \dee p_x + \dee y \wedge \dee p_y$ the Hamiltonian is
$H(x,y,p_x,p_y)$, where the equilibrium is taken to be the origin.
The Hamiltonian vector field in these coordinates is
\begin{equation} 
    \X_H = J \dee H, \quad  J = \begin{pmatrix} 0 & 1 \\ -1 & 0 \end{pmatrix}
\end{equation}
or in components
\begin{equation} 
    \dot x =  \frac{\partial H}{\partial p_x}, \,
    \dot y =  \frac{\partial H}{\partial p_y}, \quad
    \dot p_x =  -\frac{\partial H}{\partial x},\,
    \dot p_y =  -\frac{\partial H}{\partial y}\,.    
\end{equation}
Since an equilibrium of a dynamical system is a zero of the vector
field one has $\dee H(0) = 0$.  The critical value $H(0)$ of the
critical point $0$ is assumed to be $0$.  The equilibrium is
characterized by the eigenvalues $\lambda$ of its linearization $J
\dee^2 H(0)$.  If $\lambda=\alpha + \imag\omega$ is an eigenvalue, so
is $\bar \lambda$, $-\lambda$, and $-\bar\lambda$.
Real pairs of non-zero eigenvalues $\pm \alpha$ are called hyperbolic,
non-zero pure imaginary pairs $\pm \imag \omega$ are called elliptic,
and quadruples of eigenvalues with non-zero real and non-zero
imaginary part $\pm \alpha \pm \imag\omega$ are called loxodromic.
Near a loxodromic equilibrium there exist symplectic coordinates
$(\tilde x, \tilde y, \tilde p_x, \tilde p_y)$ such that the quadratic
part $H_2$ of $H$ is \cite{Williamson36,Arnold78}
\begin{equation} \label{eq:H2}
     H_2 = \alpha J_1 + \omega J_2, \qquad 
     J_1 = \tilde x \tilde p_x + \tilde y \tilde p_y, \quad
     J_2 = \tilde x \tilde p_y - \tilde y \tilde p_x \,.
\end{equation}
Introducing $z = \tilde x + \imag \tilde y$ and $p_z = \tilde p_x -
\imag \tilde p_y$ gives $\Omega = \Re(\dee z \wedge \dee p_z)$, and
$J_1 - \imag J_2 = z p_z$.  Outside $z = 0$ this extends to a
multivalued symplectic coordinate system $(\ln z, z p_z)$, and the
flow of $H_2$ is easily found to be $(z,p_z) = (\exp(\lambda t) z_0,
\exp(-\lambda t) p_{z0})$.  The quadratic Hamiltonian has a
two-dimensional stable eigenspace spanned by the eigenvectors of
eigenvalues with negative real part, and a two-dimensional unstable
eigenspace spanned by the eigenvectors of eigenvalues with positive
real part.  In each eigenspace the dynamics is that of a focus point
in the plane. The invariant manifolds of the equilibrium point of the
full Hamiltonian are tangent to these eigenspaces at the equilibrium.
In general the stable and unstable invariant manifolds will have
transverse intersections, and the system is non-integrable.

Here we are concerned with the dynamics near an equilibrium with
loxodromic eigenvalues of a Liouville integrable Hamiltonian system
with two degrees of freedom.  Under some additional hypothesis (see
below) we call such an equilibrium focus-focus point.  Let $L$ be an
independent second integral that commutes with $H$, i.e.\ the Poisson
bracket $\{ H, L \} = \Omega( J dH, J dL)$ vanishes.  The system is
called Liouville integrable if the energy-momentum map $F=(H,L)$ is
regular almost everywhere, and $\{ H, L\} = 0$.  The values in the
image of $F$ are denoted by $c = (h,l)$.  The Liouville-Arnold theorem
states that near any compact connected component of regular points
$\T_c$ of $F^{-1}(c)$ there are so called action-angle coordinates
$(\theta_1,\theta_2, I_1,I_2)$ with $\Omega = \dee \theta_1 \wedge
\dee I_1 + \dee \theta_2 \wedge \dee I_2$ and $I_1$ and $I_2$ are
commuting constants of motion with periodic flows on (and near) the
two-dimensional torus $\T_c \subseteq F^{-1}(c)$.  In these
coordinates the Hamiltonian is a function of $I_1$ and $I_2$ alone,
and the quasiperiodic flow of $H$ is the solution of $\dot \theta_i =
\partial H/\partial I_i = \omega_i(I)$, $i=1,2$.  In other words, the
Hamiltonian vectorfield $\X_H$ is a linear combination of the periodic
flows of the actions
\begin{equation} \label{eq:HofI}
      \X_H = \omega_1 \X_{I_1} + \omega_2 \X_{I_2}\,.
\end{equation}
The coefficients of the linear combination are the frequencies
$\omega_i(I)$ of the angles $\theta_i$, $i=1,2$. The frequencies
$\omega = (\omega_1,\omega_2)$ in general depend on the actions $I =
(I_1,I_2)$.  The map from actions $I$ to frequencies $\omega(I)$ is
called the frequency map.

The KAM theorem (see e.g. \cite{Arnold78,SM71}) asserts that under
small perturbations the invariant torus $\T_c$ with frequencies
$\omega(I)$ persists if 1) the frequency map is non-degenerate and 2)
the frequencies $\omega(I)$ are Diophantine.  In the isoenergetic KAM
theorem persistence of invariant tori is considered for the same
energy $h$.  Then a torus is characterized by its rotation number
$W(I)$ which is the frequency ratio $[\omega_1 : \omega_2] \in \R
P^1$.
The iso-energetic non-degeneracy condition $\partial W/\partial l\neq
0$, where the derivative is taken at constant energy, replaces the
non-degeneracy condition on the frequency map, see
e.g.~\cite{Arnold78}.  In a semi-global Poincar\'e section transversal
to $\X_H$ on $\T_c$ the iso-energetic non-degeneracy condition ensures
that the Poincar\'e map is a twist map near the invariant curve
corresponding to $\T_c$.  The twist condition $\partial W / \partial l
\not = 0$ implies that the rotation number changes between neighboring
tori.  An invariant torus $\T_c$ for which the twist vanishes is
called a twistless torus, for short.  To check the non-degeneracy of
the frequency map in an example requires tools from complex analysis,
see e.g.~Horozov's work on the spherical pendulum \cite{Horozov90}.
General results in the neighborhood of critical points of the energy
momentum map are due to Kn\"orrer \cite{Knorrer85} and Nguy\^en Ti\^en
Zung \cite{Zung96}.  We are not aware of general results about the
iso-energetic non-degeneracy condition; the spherical pendulum has
again been treated by Horozov \cite{Horozov93}.

The twist condition is only a sufficient condition for persistence,
and KAM theorems with weaker conditions exist, see \cite{Ruessmann87}.
However, the perturbed dynamics is quite unusual when the unperturbed
twistless torus is resonant, i.e.\ when it has a rational rotation
number.  A resonant torus with twist breaks into a Poincar\'e-Birkhoff
island chain, see e.g.~\cite{MH92,SM71,Birkhoff27}.  A twistless
resonant torus instead breaks into two island chains, and near the
collision of these chains interesting dynamics occurs on so called
meandering invariant curves, see \cite{HowHoh84,CGM96,Simo98,DMS98b}
and the references therein.  Our goal is to show that such dynamics
occurs near a loxodromic equilibrium of a generically perturbed
integrable system.  More precisely we will prove the following
\begin{teo} \label{teo:main}
  In every integrable Hamiltonian system with two degrees of freedom
  and a focus-focus singularity with loxodromic eigenvalues there
  exists a regular torus with vanishing twist for each value of the
  energy close to the critical one. For this energy all other tori
  close enough to the singular fibre have non-vanishing twist.
\end{teo}

This theorem might look surprising when compared to a result by Nguyen
Tien Zung \cite{Zung96}: He showed that near a focus-focus point the
Kolmogorov non-degeneracy condition is always satisfied, \emph{i.e.}
on almost all tori the frequency map is non-degenerate. Zung did not
provide a way to find the tori for which the non-degeneracy condition
fails; but, as we show in section~\ref{sec:kolmogorov}, our techniques
improve his theorem by showing that the frequency map is
non-degenerate for \emph{all} regular tori close enough to the
singular fiber.

It is well known that the Kolmogorov condition and the Twist condition
are independent, and our results show that whenever there is a
focus-focus point the twist-condition is violated while the frequency
map is non-degenerate.

\section{Rotation Number}

Let the critical value of the equilibrium be $(h,l) = (0,0)$.  For
simplicity, we assume that the corresponding singular fibre
$F^{-1}(0)$ contains only a single critical point $m$. We then need to
assume that the component of $F^{-1}(0)$ containing $m$ is compact.
Such a singularity is called of focus-focus type. Since we are
interested in tori close to the singular component of $F^{-1}(0)$ we
may disregard other connected components. Then $(0,0)$ is an isolated
critical value in the image of the energy-momentum map $F$ restricted
to those tori.
Additional structure appears near the loxodromic equilibrium by the
assumption that the system be integrable. Eliasson showed
\cite{Eliasson84} that the momenta $J_1$ and $J_2$ of the quadratic
normal form are the components of a momentum map $J = (J_1,J_2)$ of
the full system near the focus-focus point in some symplectic
coordinates.  This means that we may assume that both $H$ and $L$ (and
not only their quadratic parts) are functions of $J_1$ and $J_2$ near
$m$:
\begin{equation} \label{eq:HofJ}
   H =  \Phi \circ J \qquad \text{ and } \qquad
L = \Psi\circ J,
. 
\end{equation}
where $g=(\Phi,\Psi)$ is a local diffeomorphism of $\R^2$ near the
origin.
Recall that \eqref{eq:HofJ} is valid only near the focus-focus point
$m$.  However this gives a way of extending the momentum map $J$ to a
whole neighborhood of the singular fiber simply by enforcing
$J=g^{-1}(H,L)$.

Since $H$ is loxodromic the coefficient $\alpha$ in~\eqref{eq:H2} is
non-zero hence $\partial_1\Phi(0)\neq 0$. Then by the implicit
function theorem the additional integral $L$ can be chosen equal to
$J_2$.

Observe that $J_2$ has a $2\pi$ periodic flow near $m$, and this still
holds in a neighborhood of the singular fiber, as was shown in
\cite{Ngoc2000}. Therefore $J_2$ is a generator of an $S^1$ action on
the neighboring tori and hence is equal to one action, say $I_2$.

Instead of using the energy-momentum map $F$ it is advantageous to use
the momentum map $J$ instead, and whenever the Hamiltonian enters the
calculation to use $\Phi$.

In \cite{Ngoc2002} it is shown using only the momentum map $J$ without
the Hamiltonian (or, equivalently, with Hamiltonian equal to $J_1$)
that the $2\pi$ periodic flow $\X_{I_1}$ of the first action is given
by
\begin{equation} \label{eq:IofJ}
  2\pi \X_{I_1} = \tau_1 \X_{J_1} + \tau_2 \X_{J_2} 
\end{equation}
where the periods $\tau_i$, $i=1,2$ satisfy
\begin{equation} \label{eq:tau12}
\begin{aligned}
  \tau_1(j) =& \sigma_1(j) - \Re( \ln \zeta), \\
  \tau_2(j) =& \sigma_2(j) + \Im( \ln \zeta),
\end{aligned}
\end{equation}
and $\sigma_i(j)$, $i=1,2$, are smooth and single-valued functions
near the origin.  Here the point $j = (j_1,j_2)$ in the image of the
momentum map $J$ is identified with the complex number $\zeta = j_1 +
\imag j_2$.  The smooth contribution comes from the dynamics far away
from the focus-focus point, while the singular contribution is
obtained from the flows of $J_1$ and $J_2$ near the focus-focus point,
see \cite{Ngoc2002} for the details. The points $(\tau_1, \tau_2)$ and
$(0,2\pi)$ form a basis of the period lattice of the foliation.
Using this result we obtain
\begin{lem} \label{lem:W}
  The rotation number $W(j)$ near a focus-focus point is
   \begin{equation} \label{eq:Wjis}
      2\pi W(j) = - A(j) \Re ( \ln \zeta )  - \Im ( \ln \zeta ) + \sigma(j)
   \end{equation} 
   where $A(j)$ is smooth and determined from the Hamiltonian $H =
   \Phi \circ J$ by
   \begin{equation} 
        A(j) =  \frac{\partial_2 \Phi}{\partial_1 \Phi}(j), \quad A(0)  = \frac{\omega}{\alpha}\,,
   \end{equation} 
   and $\sigma(j)$ is a smooth (and single-valued) function near the
   origin.
\end{lem}
\proof The Hamiltonian vector field $\X_H$ is obtained from
(\ref{eq:HofJ}) as
\begin{equation} 
   \X_H = \partial_1 \Phi \X_{J_1} + \partial_2 \Phi \X_{J_2} \,.
\end{equation}
Now eliminating $\X_{J_1}$ with the aid of (\ref{eq:IofJ}) and using
$J_2 = I_2$ gives
\begin{equation} \label{eq:OmofJ}
  \X_H = \frac{2\pi}{\tau_1} \partial_1 \Phi \X_{I_1} + 
   \left( \partial_2\Phi - \frac{\tau_2}{\tau_1} \partial_1 \Phi \right)  \X_{I_2} \,.
\end{equation}
Comparing with (\ref{eq:HofI}) the frequencies $\omega_j$ can be read
off, and the rotation number as a function of $j$ is
\begin{equation} 
   W(j) =   \frac{\omega_2(j)}{\omega_1(j)} = \frac{1}{2\pi} ( \tau_1(j) A(j) - \tau_2(j) )\,.
\end{equation}
Finally using (\ref{eq:tau12}) gives the result with 
$\sigma = A \sigma_1 - \sigma_2$.  From (\ref{eq:H2})
we see that $A(0) = \omega/\alpha$.  \qed

The fact that $\Im(\ln \zeta )/2\pi$ is multivalued and increases by
$1$ upon completing a cycle $\zeta = \exp(\imag \phi)$ is a
manifestation of the fact that a simple focus-focus point has
Monodromy with index 1, see \cite{Duistermaat80,Matveev96,Zung97,Ngoc2000}.

The functions $\tau_i$ that determines the period lattice depend on the foliation alone, 
while $A(j)$ is determined by the Hamiltonian alone.  In order to understand
the behaviour of $W(j)$ near the origin we now prove the following
\begin{teo}
  \label{teo:spirals}
  There is a local diffeomorphism of $\R^2$ near the origin which is
  $C^1$ at the origin and $C^\infty$ elsewhere which maps the level
  sets of $W$ to the integral curves of the equation
  $\dot{\zeta}=-\bar{\lambda}\zeta$.  In particular these curves are
  spirals determined by the eigenvalue $\lambda$ of the focus-focus
  point.
\end{teo}

\proof We use the notation of lemma 2. The strategy is to change the
coordinates such that the level sets of $W$ in the new coordinates
satisfy a linear differential equation.  Let $\phy$ be the change of
variables defined in polar coordinates $\zeta = j_1+\imag j_2=\rho
e^{\imag\theta}$ by
\[
\left(
\begin{array}[1]{c}
\rho\\\theta
\end{array}
\right) \rightarrow 
\left(
\begin{array}[1]{l}
\tilde{\rho}=\rho^{A(j)/A(0)}\exp(-\sigma(j)/A(0))\\
\tilde{\theta}=\theta
\end{array}
\right).
\]

Writing $A(j)=A(0)+\pscal{j}{\tilde{A}(j)}$, where $\tilde{A}$ is
smooth, we have $\phy(j)=f(j)j$, where
$f:\R^2\rightarrow\R$ locally near the origin is given by
\[
f(j)=\rho^{\pscal{j}{\tilde{A}(j)}/A(0)}\exp(-\sigma(j)/A(0)).
\] 
We see that $f$ is continuous at the origin, with
$f(0)=\exp(-\sigma(0)/A(0))$, and that the partial derivatives of $f$
are of order $\ln\rho$ at the origin. Hence $\phy$ is $C^1$ at the
origin, with Jacobian
\[
{\rm d} \phy(0)= f(0) {\rm Id} = \exp(-\sigma(0)/A(0)) {\rm Id}.
\]

Under this diffeomorphism, the rotation number is simplified and
we only need to study the level sets of 
\[
2\pi W\circ\phy^{-1}(\tilde{j})=(-A(0)\Re(\ln\tilde{\zeta}) -
\Im(\ln\tilde{\zeta})).
\] 
Forget now the tildes. Viewing $W$ as a local Hamiltonian for the
standard canonical structure $\dee j_2\wedge \dee j_1$, we compute the
levels sets of $W$ as the integral curves of the associated
Hamiltonian vector field: $\frac{\dee}{\dee s}(j_1, j_2) = (
\partial_2 W, -\partial_1 W)$. Rescaling the time $s$ by
$2\pi|\zeta|^2\alpha$ and rewriting the dynamical system in complex
form gives the equation $\dot \zeta = -\bar\lambda \zeta$.

From Lemma~2 we know that $W(j)$ diverges when $j$ approaches the
origin. These statements are compatible because different branches of
the complex logarithm are involved. After each complete turn of the
spiral the rotation number jumps by one. Alternatively one could view
the rotation number as globally defined on the Riemann surface of the
complex logarithm.\qed

Recall that as coordinates in $\R^2$ we may use $j = (j_1,j_2)$ or $c
= (h,l)$, and that they are related by $l = j_2$ and $h = \Phi(j)$.
Therefore the spirals in the image of the momentum map $J$ will be
mapped to spirals in the image of the momentum map $F$.  Near the
origin the transformation between the two is the linear map $(h,l) =
(j_1 \alpha + j_2\omega, j_2)$.

\section{Vanishing Twist}

In the isoenergetic problem fixing the energy $H = h$ gives a smooth
curve $\mathscr{C}_h$ in the image of the momentum map.  The twist
condition is obtained from $W(j)$ by differentiating along this curve.
Since $\partial_1 \Phi(0) = \alpha \not = 0$ we may apply the implicit
function theorem and the curve $\mathscr{C}_h$ is the graph of a
function $j_1(j_2)$. The slope of this curve near the origin is given
by implicit differentiation of $\Phi(j_1,j_2) = h$ which yields
$\partial j_1/\partial j_2 = - \partial_2 \Phi/\partial_1 \Phi = -
A(j)$.  Therefore the twist is found to be
\begin{equation} 
   \S(j) = \left. \frac{\partial W(j)}{\partial j_2} \right|_{\mathscr{C}_h} =
    \partial_1 W  \frac{\partial j_1}{\partial j_2} + \partial_2 W
     =  -A(j)\partial_1 W + \partial_2 W\,.
\end{equation}
Expanding the remaining derivatives gives
\begin{equation} \label{eq:Sis}
  2\pi \S(j) = -A^2 \partial_1 \tau_1 + 2 A \partial_2 \tau_1 - \partial_2 \tau_2 
  -  \tau_1 A \partial_1 A  + \tau_1 \partial_2 A 
\end{equation}
where $\partial_1 \tau_2 = \partial_2 \tau_1$ has been used, which
follows from the fact that both periods can be obtained by
differentiating a single action function, see \cite{Ngoc2002}.  The
main contribution near the origin comes from differentiating the
complex logarithm in (\ref{eq:tau12}), which produces terms
proportional to $1/|j|^2$.  Therefore we now define $\tilde \S(j) =
2\pi |j|^2 \S(j)$ and have the following
\begin{lem}  \label{lem:S}
  $\tilde \S$ is a $C^1$ map at the origin that satisfies
\begin{equation} 
 \tilde \S(0) = 0, \qquad
 \partial_1 \tilde \S(0) = A^2(0) - 1, \quad
 \partial_2 \tilde \S(0) = -2 A(0)
\end{equation}
\end{lem}
\proof Using (\ref{eq:tau12}) where $\sigma_i$ is smooth we find that
the last two terms in (\ref{eq:Sis}), $ \tau_1 A \partial_1 A$ and
$\tau_1 \partial_2 A $, are both of the form $f(j) + g(j) \ln |j|$
with smooth $f$ and $g$, while the remaining first three terms are of
the form $f(j) + g(j) / |j|^2$ for some smooth $f$ and $g$.  Since
$|j|^2 \ln |j|$ is of class $C^1$ at the origin the first statement
follows.  A simple computation gives
\begin{equation} 
  \tilde \S(j) = A^2 j_1 - 2 A j_2 - j_1 + O(|j|^2 \ln |j|)\,,
\end{equation}
and the result follows.  \qed

Since the derivatives of $\tilde \S$ cannot both vanish at the origin,
$\S^{-1}(0)$ is a $C^1$ curve through the origin.  By
Lemma~\ref{lem:S} the equation for the tangent of this curve at the
origin is
\begin{equation} 
   (A(0)^2 - 1) j_1 = 2 A(0) j_2
\end{equation}
Using $h = j_1 \alpha + j_2 \omega$ at the origin to express this
tangent in $(h,j_2)$ gives
\begin{equation} \label{eq:hjis}
   h = \omega \frac{\omega^2 + \alpha^2}{\omega^2 - \alpha^2} j_2\,.
\end{equation}
When $\omega \not = 0$ the curve $\tilde \S(0)$ therefore
transversally intersects the curves $\{h={\rm const}\}$ near the
origin, and we have proven Theorem 1.  Note that the condition $\omega
\not = 0$ is only needed in the last step.  Therefore also for $\omega
= 0$ twistless tori exist near the origin, but transversality is lost.
Therefore twistless tori might exist only for positive or only for
negative values of $h$, but possibly not for both.


Different choices of actions are possible.  The semi-global $S^1$
action $L=J_2 = I_2$ is unique up to sign $\eps = \pm1$, but the other
action $I_1$ can be changed, $\tilde I_1 = \eps I_1 + k I_2$, $k \in
\Z$, so that new and old actions are related by a uni-modular
transformation.  The frequencies $\omega$ change into $\tilde \omega_1
= \eps \omega_1$ and $\tilde \omega_2 = \eps\omega_2 - k \omega_1$, so
that the rotation number transforms as $\tilde W = W - k \eps$.
Therefore a different choice of actions amounts to a different choice
of sheet of the Riemann surface of the complex logarithm.  This
argument also shows that the vanishing twist is independent of the
choice of actions.

When there are more than one critical point in the singular leave
similar results can be derived, based on modified formulas
(\ref{eq:tau12}) whose derivation is sketched in \cite{Ngoc2002}.


As an application of the above results one can consider the integrable
normal form of the Hamiltonian Hopf bifurcation, see
\cite{vanderMeer85,Duistermaat98} and the references therein.  In
\cite{DI03b} the rotation number in the compact case was computed
explicitly in terms of elliptic integrals. Expansion of these
integrals at the focus-focus point gave the above results
(\ref{eq:Wjis}) and (\ref{eq:hjis}) for the first time in this special
case.  Pictures of the spirals of constant $W$ can also be found in
\cite{DI03b}.  But the method employed could not deal with the case in
which the normal form has a non-compact singular leave. And even with
a compact singular fibre in a more complicated system this approach
might lead to hyperelliptic integrals, and their expansion at the
singular point would be quite difficult.  A prominent example of an
integrable system with a focus-focus point is the spherical pendulum,
see e.g.~\cite{CushmanBates97}.  In this case the eigenvalues are
degenerate because $\omega = 0$, and the spiral of
Theorem~\ref{teo:spirals} degenerates into a star.

\section{Kolmogorov Condition}
\label{sec:kolmogorov}

Near a focus-focus point we can improve Zung's theorem \cite{Zung96}.
\begin{teo}
  On regular tori close to the focus-focus singular fiber, the
  Jacobian determinant of the frequency map admits the following
  asymptotic expansion (recall that $\tau_1$ is of order $\ln |j|$):
  \begin{equation}
    \label{equ:det1}
    \det \frac{\partial \omega}{\partial I} = - \left( \frac{2\pi
    \alpha}{|j|\tau_1^2}\right)^2 + O(\frac{1}{|j|\tau_1^3})\,.
  \end{equation}
  In particular this shows that the Kolmogorov condition is uniformly
  satisfied for every torus near the singular fibre of a simple
  focus-focus point.
\end{teo}

\proof First notice that the Jacobian of the frequency map does not
depend on the choice of action variables and hence does not care about
monodromy; in our case, this means that it does not depend on the
determination of the complex logarithm in~\eqref{eq:tau12}.  Thus the
leading order of this Jacobian can easily be calculated. We need to
compute
\[
\frac{\partial \omega}{\partial I} =  \frac{\partial \omega}{\partial J} \frac{\partial J}{\partial I} \,.
\]
From (\ref{eq:IofJ}) we know that by definition
\[
2\pi \frac{\partial I_1}{\partial J_1} = \tau_1, \quad
2\pi \frac{\partial I_1}{\partial J_2} = \tau_2\,.
\]
Since $J_2 = I_2$ we easily find that $\det \partial I/ \partial J = \tau_1/2\pi$.
The frequencies can be read off from (\ref{eq:OmofJ}); they are
\[
\omega_1 = \frac{2\pi}{\tau_1} \Phi_1, \quad
\omega_2 = \Phi_2 - \frac{\tau_2}{\tau_1} \Phi_1\,.
\]
The subscripts of $\Phi$ denote partial derivatives. Using that $\tau_1$
is of order $\ln |j|$ and that the first derivatives of $\tau_i$ have
all leading order $1/|j|$, one easily computes:
\begin{eqnarray}
  \begin{array}{rll}
    \partial_1\omega_1 = &
    -2\pi\Phi_1\displaystyle\frac{\tau_{11}}{\tau_1^2} +
    O(1/\tau_1) & =
    \displaystyle O(\frac{1}{|j|\tau_1^2}) \\
    
    \partial_2\omega_1 = &
    -2\pi\Phi_1\displaystyle\frac{\tau_{12}}{\tau_1^2} + 
    O(1/\tau_1) & =  \displaystyle O(\frac{1}{|j|\tau_1^2}) \\
    
    \partial_1\omega_2 = & -
    \displaystyle\frac{(\tau_1\tau_{21}-\tau_2\tau_{11})}
    {\tau{}_1^2}\Phi_1 + O(1) 
    & =  \displaystyle O(\frac{1}{|j|\tau_1}) \\
    
    \partial_2\omega_2 = & -
    \displaystyle\frac{(\tau_1\tau_{22}-\tau_2\tau_{12})}
    {\tau{}_1^2}\Phi_1 + O(1) 
    & =   \displaystyle O(\frac{1}{|j|\tau_1})
  \end{array}
  \label{equ:omega}
\end{eqnarray}
Here the second index of $\tau$ denotes a partial derivative.  At
leading order $\tau_1$ can be replaced by $\ln |j|$.  The leading
order of the determinant $\det \partial \omega/\partial J$ is obtained
by taking the determinant of the matrix in which only the terms that
are given explicitly in~\eqref{equ:omega} are kept, provided one shows
\emph{a posteriori} that the result has the expected order of
$\frac{1}{|j|^2\tau_1^3} = \frac{1}{|j|^2\ln^3|j|}$.  The result is
\[
\det \frac{\partial \omega}{\partial J} = 2\pi (\tau_{11} \tau_{22} -
\tau_{12}^2) \frac{\Phi_1^2}{\tau_1^3} + O(\frac{1}{|j| \tau_1^2 } )
\,.
\]
Using the form of $\tau_1$ and $\tau_2$ we find 
\[
\det \frac{\partial \tau}{\partial J} = \tau_{11} \tau_{22} -
\tau_{12}^2 = - \frac{1}{|j|^2} + O(1/|j|)
\]
so that 
\[
\det \frac{\partial \omega}{\partial J} = - 2\pi \frac{\alpha^2}{|j|^2
  \tau_1^3} + O(\frac{1}{|j|\tau_1^2}) = -2\pi \frac{\alpha^2}{|j|^2
  \ln^3|j|} + O(\frac{1}{|j|^2\ln^4|j|})\,.
\]
Since by hypothesis $\alpha\neq 0$ the leading term has indeed the
required order.  Returning to the true actions $I$ the final result is
\[
\det \frac{\partial \omega}{\partial I} = - \left( \frac{2\pi
    \alpha}{|j|\tau_1^2}\right)^2 + O(\frac{1}{|j|\tau_1^3})\,,
\]
thereby proving the theorem. \qed

\vspace{2ex}

\noindent\emph{Remark.} In the process of finishing this paper we became
aware of a preprint by Rink~\cite{rink-focus} in which the determinant
of the Jacobian of the frequency map was calculated in a very similar
way.

\bibliographystyle{plain}

\def\cprime{$'$}

\end{document}